\documentclass{article}
\usepackage[cp1251]{inputenc}
\usepackage{amsmath}
\usepackage{amssymb}
\usepackage{amsfonts}
\usepackage{amsthm}
\usepackage{graphicx}

\def\thtext#1{
  \catcode`@=11
  \gdef\@thmcountersep{. #1}
  \catcode`@=12
}

\def\threst{
  \catcode`@=11
  \gdef\@thmcountersep{.}
  \catcode`@=12
}

 \catcode`@=11
 \def\.{.\spacefactor\@m}
 \catcode`@=12

\theoremstyle{plain}
\newtheorem{thm}{Theorem}
\newtheorem{prop}{Proposition}[section]

\newtheorem{ass}{Assertion}

\theoremstyle{definition}

\newcommand{\cP}{\mathcal{P}}
\newcommand{\cR}{\mathcal{R}}

\newcommand{\N}{\mathbb{N}}

\renewcommand{\a}{\alpha}
\renewcommand{\b}{\beta}
\newcommand{\e}{\varepsilon}
\newcommand{\s}{\sigma}

\newcommand{\gM}{\mathfrak{M}}

\newcommand{\diam}{\operatorname{diam}}
\newcommand{\dis}{\operatorname{dis}}

\renewcommand{\:}{\colon}
\renewcommand{\ss}{\subset}
\newcommand{\x}{\times}

\begin{document}
\title{The Gromov--Hausdorff Metric on the Space of Compact Metric Spaces is Strictly Intrinsic}
\author{A.~Ivanov, N.~Nikolaeva, A.~Tuzhilin\\
\small Lomonosov Moscow State University, Mechanical and Mathematical Faculty}
\date{}
\maketitle

\begin{abstract}
It is proved that the Gromov--Hausdorff metric on the space of compact metric spaces considered up to an isometry  is strictly intrinsic, i.e., the corresponding metric space is geodesic. In other words, each two points of this space  (each two compact metric spaces) can be connected by a geodesic. For finite metric spaces a geodesic is constructed explicitly.
\end{abstract}

\section{Introduction}
\markright{\thesection.~Introduction}
By  $\gM$ we denote the space of all compact metric  spaces (considered up to an isometry) endowed with Gromov--Hausdorff distance. This space is intensively investigated, and some of its properties are already known, for example, it is proved that this space is arc-connected, complete, separable, but not boundedly compact. Also, it is known that the Gromov--Hausdorff metric on this space is intrinsic\footnote{%
We have not found an exact reference for this result. It is only mentioned in Russian version of Wikipedia (but does not appear in the English or German ones). 
}. %
The aim of the present paper is to prove that this metric is also strictly intrinsic, i.e., each two points of this space can be connected by a shortest geodesic.

\section{Preliminaries}
\markright{\thesection.~Preliminaries}
All definitions and results from this Section can be found in~\cite{BurBurIva}.

Let $X$ be an arbitrary metric space. By $|xy|$ we denote the distance between its two points $x$ and $y$. For any  $x\in X$ and non-empty $A\ss X$ define $|xA|$ as the infimum of the distances $|xa|$ over all $a\in A$. For non-empty $A,\,B\ss X$ define $d(B,A)$ as the supremum of the distances $|bA|$ over all $b\in B$. At last, put $d_H(A,B)=\max\{d(A,B),d(B,A)\}$. The value $d_H(A,B)$ is referred as the \emph{Hausdorff distance}. It is well-known that $d_H$ is a metric on the set of all closed bounded subsets of the space $X$.

Let $X$ and $Y$ be metric spaces. A triplet $(X',Y',Z)$ consisting of a metric space $Z$ and two its subspaces  $X'$ and $Y'$ isometrical to $X$ and $Y$, respectively, is called a \emph{realization of the pair $(X,Y)$}. The\emph{Gromov--Hausdorff distance $d_{GH}(X,Y)$ between $X$ and $Y$} is defined as the infimum of the values  $r$ such that there exists a realization $(X',Y',Z)$ of the pair $(X,Y)$ such that $d_H(X',Y')\le r$. It is well-known that the function $d_{GH}$ is a metric on the set $\gM$, and that this metric is intrinsic.

Recall that a binary \emph{relation\/} on the sets $X$ and $Y$ is an arbitrary subset of the direct product $X\x Y$. By $\cP(X,Y)$ we denote the set of all relations on $X$ and $Y$. Let $\pi_X\:X\x Y\to X$ and $\pi_Y\:X\x Y\to Y$ be the canonical projections, i.e., $\pi_X(x,y)=x$ and $\pi_Y(x,y)=y$. In the same way we denote the restrictions of the canonical projections to each relation $\s\ss\cP(X,Y)$.

A relation $R$ on  $X$ and $Y$ is called  a \emph{correspondence\/} or \emph{left-right-total}, if the restrictions of the canonical projections  $\pi_X$ and $\pi_Y$ to $R$ are surjective. In other words, for any $x\in X$ there exists  $y\in Y$ such that  $x$ and $y$ are in relation $R$ and, conversely, for any $y\in Y$ there exists $x\in X$ such that $x$ and $y$ are in relation  $R$. By $\cR(X,Y)$ we denote the set of all correspondences on $X$ and $Y$.

Let $X$ and $Y$ be metric spaces, then for each relation $\s\in\cP(X,Y)$ its \emph{distortion\ $\dis\s$} is defined as follows:
$$
\dis\s=\sup\Bigl\{\bigl||xx'|-|yy'|\bigr|: (x,y)\in\s,\ (x',y')\in\s\Bigr\}.
$$
If $f\:X\to Y$ is a mapping, then its \emph{distorsion $\dis f$} is defined as the distortion of its graph $\bigl\{(x,f(x))\mid x\in X\Bigr\}\subset X\x Y$.

The next result is well-known (Theorem~7.3.25 in~\cite{BurBurIva}).

\begin{prop}\label{th:GH-metri-and-relations}
For any metric spaces $X$ and $Y$ the equality
$$
d_{GH}(X,Y)=\frac12\inf\bigl\{\dis R\mid R\in\cR(X,Y)\bigr\}
$$
is valid.
\end{prop}

If $X$ and $Y$ are finite metric spaces, then the set $\cR(X,Y)$ is finite, and hence, there exists $R\in\cR(X,Y)$ such that $d_{GH}(X,Y)=\frac12\dis R$. Each such $R$ is referred as \emph{optimal}.

We also need some other auxiliary results. Recall that a subset $S$ of a metric space $X$ is called an \emph{$\e$-net\/} in $X$, if for each $x\in X$ there exists some $s\in S$ with $|sx|<\e$. 

\begin{prop}\label{prop:net}
Let $X$ be an arbitrary metric space, $Y$ be its non-empty subset $X$, and $S$ be some $\e$-net in $X$. Then there exists a $(2\e)$-net in $Y$, such that its cardinality is less than or equal to the cardinality of  the set $S$.
\end{prop}

Let $M\ss\gM$ be a family of compact metric spaces. It is said to be \emph{uniformly totally bounded}, if the following two conditions hold:
\begin{enumerate}
\item there exists a number $D\ge0$ such that for any $X\in M$ the inequality $\diam X\le D$ is valid (i.e., the diameters of the spaces from $M$ are uniformly bounded);
\item for any $\e>0$ there exists $N(\e)\in\N$ such that each $X\in M$ contains an $\e$-net consisting of at most  $N(\e)$ points (i.e., the sizes of the $\e$-nets of the spaces from $M$ are uniformly bounded).
\end{enumerate}

The following criterion is well-known (Theorem~7.4.15 in~\cite{BurBurIva}).

\begin{prop}\label{prop:precomp}
A family $M\ss\gM$ is precompact, if and only if it is uniformly totally bounded.
\end{prop}

\section{A Geodesic Connecting Finite Metric Spaces}
\markright{\thesection.~A Geodesic Connecting Finite Metric Spaces}

Let $X$ and $Y$ be two finite metric spaces, and $R$ be an optimal correspondence on them. Then $d_{GH}(X,Y)=\frac12\dis R$, due to definitions.

Fix an arbitrary $0\le\a\le 1$, and let $\b=1-\a$. We put
$$
\bigl|(x,y)(x',y')\bigr|_\a=\a|xx'|+\b|yy'|.
$$
Evidently, the resulting function is positively defined  and symmetric. Also, the triangle inequality is valid:
\begin{multline*}
\bigl|(x,y)(x',y')\bigr|_\a+\bigl|(x',y')(x'',y'')\bigr|_\a=\a|xx'|+\b|yy'|+\a|x'x''|+\b|y'y''|\ge\\ \ge \a|xx''|+\b|yy''|=\bigl|(x,y)(x'',y'')\bigr|_\a.
\end{multline*}
Thus, $|\cdot|_\a$ is a metric on $R$ for each $0<\a<1$, and $|\cdot|_\a$ is a semi-metric for $\a=0,\,1$. Notice that identifying the points on zero distance for the cases $\a=0,\,1$ we obtain metric spaces isometric to $X$ and $Y$, respectively. 

By $R_x\ss X\x R$ and $R_y\ss R\x Y$ we denote the correspondences defined as follows:
$$
R_x=\Bigl\{\bigl(x,(x,y)\bigr):x\in X,\,(x,y)\in R\Bigr\},\ \ R_x=\Bigl\{\bigl((x,y),y\bigr):(x,y)\in R,\,y\in Y\Bigr\}.
$$
Calculate the distortions of those correspondences. We obtain:
\begin{multline*}
\dis R_x=\max\biggl\{\Bigl||xx'|-\a|xx'|-\b|yy'|\Bigr|:x,\,x'\in X,\ (x,y),\,(x',y')\in R\biggr\}=\\ =
\b\max\Bigl\{\bigl||xx'|-|yy'|\bigr|: (x,y),\,(x',y')\in R\Bigr\}=\b\dis R.
\end{multline*}
Similarly, $\dis R_y=\a\dis R$.

Therefore, $d_{GH}(X,R)\le\b\, d_{GH}(X,Y)$ and $d_{GH}(R,Y)\le\a\, d_{GH}(X,Y)$, and hence, $d_{GH}(X,R)=\b\, d_{GH}(X,Y)$, and $d_{GH}(R,Y)=\a\, d_{GH}(X,Y)$, i.e., $\bigl(R,|\cdot|_\a\bigr)$ lies between $X$ and $Y$ for any $\a\in[0,1]$. In particular, the space $\bigl(R,|\cdot|_{1/2}\bigr)$ is a midpoint between $X$ and $Y$.

Thus, the following result is proved.

\begin{ass}
For any two finite metric spaces $X$ and $Y$ there exists a midpoint between $X$ and $Y$ in $\gM$. Moreover, if $d=d_{GH}(X,Y)$, then the mapping $g\:[0,d]\to\gM$ defined as $g(t)=\bigl(R,|\cdot|_{t/d}\bigr)$ is a geodesic connecting $X$ and $Y$.
\end{ass}

The midpoint between finite metric spaces constructed above by an optimal correspondence $R$ is referred as  the \emph{canonical midpoint constructed by the optimal correspondence $R$\/} and is denoted by the same symbol  $R$.

\section{Proof of the Main Theorem}
\markright{\thesection.~Proof of the Main Theorem}

Now let $X$ and $Y$ be arbitrary compact metric spaces. For each $n\in\N$ by  $X_n$ and $Y_n$ we denote some finite $1/n$-nets in $X$ and $Y$, respectively, and let $R_n$ be a canonical midpoint between $X_n$ and $Y_n$. Show that the set $\{R_n\}_{n\in\N}\ss\gM$ is precompact. Since $X_n\to X$ and $Y_n\to Y$, then, due to Proposition~\ref{prop:precomp}, there exists a number $D$ such that $\diam X_n\le D$ and $\diam Y_n\le D$ for all $n\in\N$, and also for each $\e>0$ there exists $N(\e)>0$ such that each $X_n$ and $Y_n$ contain $\e$-nets $X'_n$ and $Y'_n$ consisting of at most $N(\e)$ points.

Define a distance function on $X_n\x Y_n$ using the same idea as in the case of the correspondences $R_n$, namely, put
$$
\bigl|(x,y),(x',y')\bigr|=\frac12\bigl(|xx'|+|yy'|\bigr).
$$
It is clear that the restriction of this distance function onto $R_n$ coincides with the distance defined on $R_n$ above as on the canonical midpoint. Besides, $\diam R_n\le\diam(X_n\x Y_n)\le D$, hence the set $\{R_n\}$ is uniformly bounded. At last, $X'_n\x Y'_n$ is an $\e$-net in $X_n\x Y_n$ consisting of at most $N(\e)^2$ points. Therefore $R_n$ contains a $(2\e)$-net consisting of at most $N(\e)^2$ points, due to Proposition~\ref{prop:net}. Thus, the conditions of the Gromov precompact criterion (Proposition~\ref{prop:precomp}) are valid, hence, the set $\{R_n\}$ is precompact.

Without loss of generality, assume that the sequence $R_n$ converges in $\gM$ to some compact metric space $R$. Continuity of the distance function implies that $R$ is a midpoint between $X$ and $Y$. Thus, the following result is proved.

\begin{ass}
For any two compact metric spaces, there exists a midpoint between them in $\gM$.
\end{ass}

The following result is well-known (see Theorem 2.4.16 in~\cite{BurBurIva}).

\begin{prop}\label{prop:mid}
Let $X$ be a complete metric space. If for any two points $a$ and $b$ from $X$ there exists a midpoint between $a$ and $b$, then $X$ is a geodesic metric space.
\end{prop}

The completeness of $\gM$ together with Proposition~\ref{prop:mid} implies the Main Theorem.

\begin{thm}
The Gromov--Hausdorff metric on the space $\gM$ of compact metric spaces considered up to an isometry is stricly intrinsic, i.e., the space $\gM$ is geodesic.
\end{thm}


\begin{thebibliography}{99}
\bibitem{BurBurIva} 
D.~Burago, Yu.~Burago, and S.~Ivanov, {\it A Course in Metric Geometry}, Graduate Studies in Mathematics, vol~33 (A.M.S., Providence, RI, 2001; Russian Edition: IKI, Moscow, Izhevsk, 2004).

\end{thebibliography}
\end{document}